\documentclass[reqno]{amsart}
\usepackage{mathrsfs,amssymb}                   % to get fancy script fonts
\numberwithin{equation}{section}

\theoremstyle{plain}
\newtheorem{thm}{Theorem}[section]
\newtheorem{prop}[thm]{Proposition}
\newtheorem{lemma}[thm]{Lemma}
\newtheorem{cor}[thm]{Corollary}

\theoremstyle{definition}
\newtheorem{defn}[thm]{Definition}
\newtheorem{example}[thm]{Example}

\newcommand{\Case}[2]{\noindent \textit{Case~#1:~#2}}

\newcommand{\numd}[2]{[#1]_{#2}}        % number of elements in #1 at depth #2
\newcommand{\otherwise}{\text{otherwise}}    % for cases
\newcommand{\Hadmat}[2]{H^{#1,#2}}      % Hadamard matrix
\newcommand{\Id}{{\bf 1}}               % identity of M_{12}
\newcommand{\Cayley}{\Gamma}

\newcommand{\ptuple}{{\bf p}}
\newcommand{\qtuple}{{\bf q}}
\newcommand{\Sym}{{\mathfrak S}}        % symmetric group
\newcommand{\Alt}{{\mathfrak A}}        % alternating group
\newcommand{\isom}{\cong}
\newcommand{\flip}[1]{\underline{#1}}
\newcommand{\comment}[1]{}
\newcommand{\defterm}[1]{{\em #1}}
\newcommand{\Prp}{{\mathbb P}_3}        % projective plane of order 3
\newcommand{\Ff}{{\mathbb F}}           % field of 3 elements
\newcommand{\Zz}{{\mathbb Z}}
\newcommand{\GM}{G_{\rm bas}}           % simple game group
\newcommand{\GS}{G_{\rm sgn}}           % signed game group
\newcommand{\GD}{G_{\rm dual}}          % signed game group
\newcommand{\0}{^{\phantom 0}}          % for spacing
\newcommand{\C}{{\mathcal C}}           % proto-Golay code
\newcommand{\G}{{\mathcal G}}           % Golay code
\newcommand{\LL}{{\mathcal L}}          % lines of \Prp
\newcommand{\PP}{{\mathcal P}}          % points of \Prp
\newcommand{\sm}{\setminus}
\newcommand{\st}{~\mid~} % "such that"
\newcommand{\x}{\times}

\newcommand{\pline}[1] {\overline{#1}}         % use for numbered lines
\newcommand{\twopt}[2] {\overline{#1#2}}        % the line defined by two points
\newcommand{\starpt}[1] {\LL(#1)}               % set of lines containing a given point

\newcommand{\Aut}{\mathop{\rm Aut}\nolimits}
\newcommand{\PGL}{\mathop{\rm PGL}\nolimits}
\newcommand{\Inn}{\mathop{\rm Inn}\nolimits}
\newcommand{\Supp}{\mathop{\rm Supp}\nolimits}
\newcommand{\minwt}{\wt_{\min}}
\newcommand{\wt}{\mathop{\rm wt}\nolimits}
\newcommand{\pmv}[2]{[#1,#2]}            % "point move"
\newcommand{\lmv}[2]{[#1,#2]}            % "line move"

\begin{document}
\title{The Mathieu group $M_{12}$ and its pseudogroup extension $M_{13}$}

\author{John H.\ Conway}
\address{Department of Mathematics\\
Princeton University\\
Princeton, NJ 08544}

\author{Noam D.\ Elkies}
\address{Department of Mathematics\\
Harvard University\\
Cambridge, MA 02138}
\email{elkies@math.harvard.edu}

\author{Jeremy L.\ Martin}  
\address{Department of Mathematics\\
University of Kansas\\
Lawrence, KS 66045}
\email{jmartin@math.ku.edu}

\keywords{Mathieu group, finite projective plane, Golay code, Hadamard matrix}
\subjclass[2000]{Primary 20B25; Secondary 05B25, 51E20, 20B20}
\thanks{Third author supported in part by an NSF Postdoctoral Fellowship}

\begin{abstract}
We study a construction of the Mathieu group $M_{12}$
using a game reminiscent of Loyd's ``15-puzzle''.
The elements of $M_{12}$ are realized as permutations
on~$12$ of the~$13$ points of the finite projective plane of order~$3$.
There is a natural extension to a ``pseudogroup'' $M_{13}$
acting on all~$13$ points, which exhibits a limited form
of sextuple transitivity.  Another corollary of the construction
is a metric, akin to that induced by a Cayley graph,
on both $M_{12}$ and $M_{13}$.  We develop these results,
and extend them to the double covers and automorphism groups
of $M_{12}$ and $M_{13}$, using the ternary Golay code and
$12 \x 12$ Hadamard matrices.
In addition, we use experimental data on the quasi-Cayley metric
to gain some insight into the structure
of these groups and pseudogroups.
\end{abstract}

\maketitle

%%%%%%%%%%%%%%%%%%%%%%%%%%%%%%%%%%%%%%%%%%%%%%%%%%%%%%%%%%%%%%%%%%%%%%%%
%%%%%%%%%%%%%%%%%%%%%%%%%%%%%%%%%%%%%%%%%%%%%%%%%%%%%%%%%%%%%%%%%%%%%%%%
%%%%%%%%%%%%%%%%%%%%%%%%%%%%%%%%%%%%%%%%%%%%%%%%%%%%%%%%%%%%%%%%%%%%%%%%

\section{Introduction} \label{intro-section}

Sam Loyd's classic \defterm{15-puzzle} consists of~$15$
numbered tiles placed in a $4 \x 4$ square grid, with one square,
the \defterm{hole}, left empty.  To solve the puzzle, one slides
the tiles around the grid until they are in a specified order.
Each sequence of slides induces a permutation
in the symmetric group $\Sym_{16}$.
The permutations arising from \defterm{closed} sequences
of slides---that is, sequences that return the hole to its initial
location---form a subgroup of the symmetric group $\Sym_{15}$,
with the group operation given by concatenation of sequences.
This subgroup is known to be the alternating group $\Alt_{15}$
(see~\cite{{Archer:1999}}).

We study an analogous game, first mentioned in~\cite{Conway:1987},
in which the $4 \x 4$ grid of Loyd's puzzle is replaced by $\Prp$,
the projective plane of order~$3$.  In the ``basic game'',
we place numbered counters on~$12$ of the~$13$ points of $\Prp$,
leaving a hole at the thirteenth point.
The elementary move, analogous to sliding an adjacent tile
to the empty square in Loyd's puzzle, is a double transposition
taking place in a line containing the hole.  
The \defterm{basic $\Prp$-game group} $\GM$
consists of the permutations of the~$12$ counters
coming from closed move sequences.
We shall prove that $\GM$ is isomorphic to the Mathieu group $M_{12}$.

We give the name $M_{13}$ to the set of permutations
induced by arbitrary (not necessarily closed) move sequences.
This is a subset of $\Sym_{13}$, but is not a group,
because concatenation of arbitrary move sequences
is not always allowed.  Specifically, a move sequence
moving the hole from $p$ to~$q$
may be followed by one taking the hole from $r$ to~$s$
if and only if $q=r$.

The $\Prp$-game can be extended in two ways.
First, we can make the counters two-sided
and modify the definition of a move to flip certain counters.
We study this ``signed game'' in Section~\ref{signed-section}.
The group~$\GS$ resulting from this change
is the nontrivial double cover $2M_{12}$ of the Mathieu group
(see~\cite[pp.~31--32]{Conway:1985}), realized as the automorphism group
of the ternary Golay code ${\mathscr C}_{12}$.
The set $2M_{13}$ of all reachable signed permutations
 is thus a double cover of $M_{13}$.

A second way to extend the basic game
is to place a second set of counters on the lines of $\Prp$.
We study this ``dualized game'' in Section~\ref{dualized-section}.
This approach yields another proof (using Hadamard matrices)
that the group $\GM$ is isomorphic to $M_{12}$;
in addition, we obtain a concrete interpretation
of an outer automorphism of $M_{12}$.

$M_{12}$ is unique among groups
in having a faithful and sharply quintuply transitive action
on a $12$-element set.  Our construction of $M_{13}$ suggests
the following question: does $M_{13}$ have a sextuply transitive
``action'' on the $13$-element set $\Prp$?
In general, the answer is no, but $M_{13}$ does exhibit
some limited forms of sextuple transitivity, which we describe
in Section~\ref{sextuple-section}.

In Section~\ref{metric-section},
we study the quasi-Cayley metric on $M_{12}$ and its extensions,
in which the distance $d(\sigma,\tau)$ between two permutations
is the minimum number of moves of the basic or signed game
needed to realize $\sigma^{-1}\tau$.  We programmed a computer
to generate lists of all positions of each of the versions of the
$\Prp$-game.  The data in these lists provides a starting point
for investigation of various aspects of the structure
of the groups and pseudogroups.
For instance, the $9$-element \defterm{tetracode}
(see, e.g., ~\cite[p.~81]{Conway:1999}) appears as the subgroup of $M_{12}$
consisting of the starting position of the basic game,
together with the $8$ positions at maximal distance from it.

The construction of $M_{13}$ was first given by the first
author in \cite{Conway:1987}.
Much of the material of this paper comes from the third author's
undergraduate thesis~\cite{Martin:1996}, written under the direction
of the second author.

%%%%%%%%%%%%%%%%%%%%%%%%%%%%%%%%%%%%%%%%%%%%%%%%%%%%%%%%%%%%%%%%%%%%%%%%
%%%%%%%%%%%%%%%%%%%%%%%%%%%%%%%%%%%%%%%%%%%%%%%%%%%%%%%%%%%%%%%%%%%%%%%%
%%%%%%%%%%%%%%%%%%%%%%%%%%%%%%%%%%%%%%%%%%%%%%%%%%%%%%%%%%%%%%%%%%%%%%%%

\section{The basic and signed $\Prp$-games} \label{game-section}
\subsection{Finite projective planes} \label{pthree-subsection}

We begin by reviewing the definitions and facts we will need
concerning~$\Prp$.

Let $\PP = \{p,q,\dots\}$ be a finite set of points
and $\LL=\{\ell,m,\dots\}$ be a finite set of lines.
Each line may be regarded as a set of points;
we write $\starpt{p}$ for the set of lines containing a point~$p$.

\begin{defn} \label{define-P3}
Let $n \geq 2$.  The pair $(\PP,\LL)$
is a \defterm{projective plane of order~$n$}
if the following conditions hold:
        \begin{enumerate}
        \item $|\PP| = |\LL| = n^2+n+1$.
        \item $|\starpt{p}| = |\ell| = n+1$ for every $p \in \PP$ and $\ell 
                \in \LL$.
        \item Two distinct points $p,q$ determine a unique line $\twopt{p}{q}$,
                that is, $\starpt{p} \cap \starpt{q} = \{\twopt{p}{q}\}$.
        \item Two distinct lines $\ell,m$ determine a unique point, that is,
                $|\ell \cap m| = 1$.  (We often write $\ell \cap m$ for the
                unique intersection point.)
        \end{enumerate}
\end{defn}

The most familiar example of a finite projective plane
is the \defterm{algebraic projective plane ${\mathbb P}_n$},
where $n$ is a prime power.  ${\mathbb P}_n$ is defined as follows.
Let $\Ff_n$ be the field of $n$ elements,
and let $\PP$ and $\LL$ be respectively the sets of
$1$- and $2$-dimensional vector subspaces of $(\Ff_n)^3$,
with incidence given by inclusion.  
Then ${\mathbb P}_n = (\PP,\LL)$ is a projective plane of order~$n$.

An \defterm{oval} (resp.\ \defterm{hyperoval}\/)
in a projective plane of order~$n$
is a set of $n+1$ (resp.~$n+2$) points, no three of which are collinear.
For example, a smooth conic in an algebraic projective plane is an oval.
Also, there is no such thing as a ``hyperhyperoval'',
for the following reason.  Let $S$\/ be a set of $n+3$ or more points
in a projective plane of order~$n$, and $p \in S$.
There are $n+1$ lines containing $p$, so by the pigeonhole principle
there exist two distinct points in $S \sm \{p\}$
which are collinear with $p$.

For the rest of the paper, we shall be exclusively concerned
with the algebraic projective plane~$\Prp$.
Note that $\Prp$ has~$13$ points, each lying on~$4$ lines,
and~$13$ lines, each containing~$4$ points.
By elementary counting, $\Prp$ has
$13 \cdot 12 \cdot 9 \cdot 4 = 5616$ ordered ovals.

\begin{prop} \label{P3-uniqueness}
\begin{enumerate}
\item Up to isomorphism, $\Prp$ is the unique finite projective
plane of order~$3$.
\item The automorphism group $\Aut(\Prp)$
  acts sharply transitively on ordered ovals; in particular,
  $|\Aut(\Prp)| = 13 \cdot 12 \cdot 9 \cdot 4 = 5616$.
\item $\Prp$ contains no hyperovals.
\end{enumerate}
\end{prop}

These facts are well known (see, e.g.,~\cite{Cameron:1991}).
We prove them here by constructing an explicit labelling
for the points and lines of $\Prp$, which we will use in the proofs
of Proposition~\ref{omnibus-prop} and Theorem~\ref{UR-thm}.

\begin{proof}
Let $\Prp$ be a projective plane of order~$3$.
Let $O=(q_1,q_2,q_3,q_4)$ be an ordered oval in~$\Prp$.
We will show that each point and each line of $\Prp$ is determined uniquely
as a function of $q_1,q_2,q_3,q_4$.  This will show that if $\Prp'$
is any other projective plane of order~$3$,
and $(q'_1,q'_2,q'_3,q'_4)$ is an ordered oval in $\Prp'$,
then there is a unique isomorphism from $\Prp$ to $\Prp'$
sending $(q_1,q_2,q_3,q_4)$ to $(q'_1,q'_2,q'_3,q'_4)$.
In particular,
it will follow that all projective planes of order~$3$ are isomorphic,
with an isomorphism group of order $13 \cdot 12 \cdot 9 \cdot 4$. 
We already know one, namely the algebraic projective plane over $\Ff_3$;
and we know that its automorphism group contains $\PGL_3(\Ff_3)$,
a group of order
$(26 \cdot 24 \cdot 18)/2 = 13 \cdot 12 \cdot 9 \cdot 4$.
Hence $\Prp$ is the unique projective plane of order~$3$,
and $\PGL_3(\Ff_3)$ is its full automorphism group.
Along the way, we will show that $\Prp$ has no hyperovals.

The points of $\Prp$ include:
\begin{itemize}
\item the four points $q_1,q_2,q_3,q_4$;
\item the three points $r_1 = \twopt{q_1}{q_2} \cap \twopt{q_3}{q_4}$,
  $r_2 = \twopt{q_1}{q_3} \cap \twopt{q_2}{q_4}$,
  $r_3 = \twopt{q_1}{q_4} \cap \twopt{q_2}{q_3}$; and
\item the six points $s_{ij}$ ($1 \leq i < j \leq 4$)
   lying on just one of the secant lines $\twopt{q_i}{q_j}$ to $O$.
\end{itemize}
(Recall that a \defterm{secant} to an oval
is a line determined by two of its points;
thus each secant $\twopt{q_i}{q_j}$ to $O$ contains
 two $q_i$'s, one $r_i$, and one other point, which we call $s_{ij}$.)
This accounts for all $4+3+6=13$ points,
and in particular shows that there are no hyperovals.

We have accounted for six lines, namely the secants to $O$.
There are also the four \defterm{tangents} to $O$,
the lines which pass through exactly one of its points.
The tangent to $O$ at $q_1$ intersects each of the lines
$\twopt{q_2}{q_3}$, $\twopt{q_2}{q_4}$, $\twopt{q_3}{q_4}$
in a point not on any other secant, which is that line's fourth point.
We have thus identified the~$4$ points on each tangent.
Three additional lines remain to be identified.

We claim that $r_1,r_2,r_3$ are not collinear,
and thus that the lines through pairs in $r_1,r_2,r_3$
complete the roster of lines of $\Prp$.
Consider for instance $s_{12}$, the fourth point on $\twopt{q_1}{q_2}$.
It lies also on the tangents at $q_3$ and $q_4$.
The points on these three lines are:
$s_{12}$ itself; $q_1,q_2,r_1$;
$q_3,s_{14},s_{24}$; $q_4,s_{13},s_{23}$.
Hence the remaining line through $s_{12}$
goes through $r_2$, $r_3$, and $s_{34}$.
This identifies the line $\twopt{r_2}{r_3}$
and shows that it does not contain $r_1$.
Likewise, we find that the line $\twopt{r_1}{r_2}$
contains $s_{14}$ and $s_{23}$,
while $\twopt{r_1}{r_3}$ contains $s_{13}$ and $s_{24}$.

We have now identified all $13$ points and all $13$ lines of $\Prp$
and their incidence relation, as desired.  To summarize, the lines are:
\begin{equation} \label{oval-labelling}
\begin{array}{lll}
 \{q_1,\,q_2,\,r_1,\,s_{12}\},\, & \{q_1,\,s_{23},\,s_{24},\,s_{34}\},\,
    & \{r_1,\,r_2,\,s_{14},\,s_{23}\}, \\
 \{q_1,\,q_3,\,r_2,\,s_{13}\},\, & \{q_2,\,s_{13},\,s_{14},\,s_{34}\},\,
   & \{r_1,\,r_3,\,s_{13},\,s_{24}\}, \\
 \{q_1,\,q_4,\,r_3,\,s_{14}\},\, & \{q_3,\,s_{12},\,s_{14},\,s_{24}\},\,
   & \{r_2,\,r_3,\,s_{12},\,s_{34}\}, \\
 \{q_2,\,q_3,\,r_3,\,s_{23}\},\, & \{q_4,\,s_{12},\,s_{13},\,s_{23}\}, \\
 \{q_2,\,q_4,\,r_2,\,s_{24}\}, \\
 \{q_3,\,q_4,\,r_1,\,s_{34}\}. \\
\end{array}
\end{equation}
%\begin{equation} \label{oval-labelling}
%\begin{array}{lll}
% \{q_1,\,q_2,\,r_1,\,s_{12}\},\, & \{q_1,\,s_{23},\,s_{24},\,s_{34}\},\,
%    & \{r_1,\,r_2,\,s_{12},\,s_{34}\}, \\
% \{q_1,\,q_3,\,r_2,\,s_{13}\},\, & \{q_2,\,s_{13},\,s_{14},\,s_{34}\},\,
%   & \{r_1,\,r_3,\,s_{13},\,s_{24}\}, \\
% \{q_1,\,q_4,\,r_3,\,s_{14}\},\, & \{q_3,\,s_{12},\,s_{14},\,s_{24}\},\,
%   & \{r_2,\,r_3,\,s_{14},\,s_{23}\}, \\
% \{q_2,\,q_3,\,r_3,\,s_{23}\},\, & \{q_4,\,s_{12},\,s_{13},\,s_{23}\}, \\
% \{q_2,\,q_4,\,r_2,\,s_{24}\}, \\
% \{q_3,\,q_4,\,r_1,\,s_{34}\}. \\
%\end{array}
%\end{equation}
\end{proof}

The definition of $\Prp$ is \defterm{self-dual}
in the sense that interchanging the terms ``point'' and ``line''
preserves the definition.
One can label the points by $\{0,1,\dots,12\}$
and the lines by $\{\pline{0},\pline{1},\dots,\pline{12}\}$
such that the lines containing the point~$x$
have the same labels as the points of the line~$\pline{x}$.  
For future reference, we give one such labelling:
\begin{equation} \label{dual-labelling}
  \begin{array}{llll}
  \pline{0} = \{0,1,2,3\}, &
  \pline{1} = \{0,4,5,6\}, &
  \pline{2} = \{0,9,10,11\}, &
  \pline{3} = \{0,7,8,12\}, \\
  \pline{4} = \{1,4,8,9\}, &
  \pline{5} = \{1,6,7,11\}, &
  \pline{6} = \{1,5,10,12\}, &
  \pline{7} = \{3,5,8,11\}, \\
  \pline{8} = \{3,4,7,10\}, &
  \pline{9} = \{2,4,11,12\}, &
  \pline{10} = \{2,6,8,10\}, &
  \pline{11} = \{2,5,7,9\}, \\
  \pline{12} = \{3,6,9,12\}.
  \end{array}
\end{equation}

%%%%%%%%%%%%%%%%%%%%%%%%%%%%%%%%%%%%%%%%%%%%%%%%%%%%%%%%%%

\subsection{The basic $\Prp$-game and $M_{13}$}
\label{basic-game-subsection}

We now describe a ``game'' similar to Loyd's $15$-puzzle,
but played on the projective plane $\Prp$ rather than a square grid.
Throughout, we use the self-dual labelling~(\ref{dual-labelling}).

To start the game, we place counters numbered $1,\dots,12$
on the respective points of $\Prp$, leaving a hole at the point~$0$.
A move of the game is defined as follows.  
Suppose that the hole is a point~$p$
and that $\ell = \{p,q,r,s\}$ is a line of $\Prp$.
Then the move $\pmv{p}{q}$
consists of moving the counter on~$q$ to~$p$
and interchanging the counters on $r$ and $s$.
This notation is justified because the pair $\{r,s\}$
is uniquely determined by the points $p$ and $q$,
by the definition of a projective plane of order~$3$.
Moreover, the move $\pmv{p}{q}$ transfers the hole from~$p$ to~$q$,
so the next move must be of the form $\pmv{q}{t}$ for some $t$.
In general, a sequence of moves can be given by specifying
the path traversed by the hole:
\begin{equation} \label{move-sequence}
  [p_0,\,p_1,\,\dots,\,p_n] \ = \ \pmv{p_{n-1}}{p_n} \,\circ\, \dots \,\circ\,
  \pmv{p_1}{p_0}.
\end{equation}
By convention, the move $\pmv{p}{p}$ is trivial, so there are~$12$ 
nontrivial legal moves playable from each position of the game.

The move $\pmv{p}{q}$ may be regarded as inducing the permutation
$(p~q)(r~s) \in \Sym_{13}$,
and a move sequence such as that of~(\ref{move-sequence})
induces the permutation
$$
(p_{n-1}~p_n)(q_n~r_n) \cdots (p_0~p_1)(q_1~r_1),
$$
where $q_i,r_i$ are the other two points
on the line $\twopt{p_{i-1}}{p_i}$ for each $i$ 
(assuming that the sequence contains no trivial moves).
Here multiplication proceeds right to left,
as is usual for permutations.

\begin{example} \label{basic-example}
Consider the path $[0,6,12,1,8,0]$.
Since the points $0$ and $6$ are collinear with $4$ and $5$,
the first move $\pmv{0}{6}$ induces the permutation
$(0~6)(4~5)$.  The permutation induced by the entire path is
\begin{multline*}
(0~8)(7~12) \cdot (1~8)(4~9) \cdot (1~12)(5~10) \cdot (6~12)(3~9) 
  \cdot (0~6)(4~5) \\
=~ (1~7~12~6~8)(3~4~10~5~9). 
\end{multline*}
\end{example}

Two paths are called \defterm{equivalent}
if they induce the same permutation.
We readily check that if $p,q,r$ are collinear
then the paths $[p,q,r]$ and $[p,r]$ are equivalent.
It follows that every path
is equivalent to a path of equal or shorter length
in which no three consecutive points are collinear;
we say that such a path is \defterm{nondegenerate}.

We say that a path $[p_0,\dots,p_n]$ is \defterm{closed}
  if $p_0=p_n$.
The set of permutations induced by closed move sequences
with $p_0=p_n=0$ is a subgroup of the symmetric group
$\Sym_{\PP \sm \{0\}} = \Sym_{12}$.
We call this subgroup the \defterm{basic $\Prp$-game group}~$\GM$,
and denote its identity element by $\Id$.

The permutations realized by move sequences
taking the hole from~$p$ to~$q$
constitute a double coset of $\GM$ in $\Sym_{\PP}$,
namely $\pmv{0}{q}\:\GM\:\pmv{p}{0}$.
In the case that $p=q$, this double coset is a group
which we call the \defterm{$q$-conjugate of~$\GM$}.

We denote by $M_{13}$ the set of all
(not necessarily closed) move sequences with $p_0=0$.
This name will be justified when we prove that
$\GM$ is isomorphic to the Mathieu group $M_{12}$.
Note that $M_{13}$ is not a group:
the moves available in a given position
depend on the location of the hole,
so concatenation of move sequences is not always allowed.
Rather, $M_{13}$ is a disjoint union of cosets
of~$\GM$ in \hbox{$\Sym_\PP = \Sym_{13}$}.

%%%%%%%%%%%%%%%%%%%%%%%%%%%%%%%%%%%%%%%%%%%%%%%%%%%%%%%%%%

\subsection{The signed $\Prp$-game} \label{signed-game-section}

We now describe the \defterm{signed $\Prp$-game},
an extension of the $\Prp$-game
in which each counter has two distinguishable sides.
Suppose that the hole is at $p \in \PP$
and that $\ell = \{p,q,r,s\} \in \LL$.
The move $\pmv{p}{q}$ of the signed game moves the counter on~$q$ to~$p$
and interchanges the counters on~$r$ and~$s$,
but it also flips over the counters on~$r$ and~$s$.  
Now a move sequence may be regarded
as inducing a signed permutation on~$\PP$
(that is, an element of the wreath product $\Zz/2\Zz \ \wr \ \Sym_{\PP}$).

\begin{example} \label{signed-example}
The path $[0,6,12,1,8,0]$ induces the permutation
\begin{multline*}
(0~8)(\flip{7}~\flip{12}) \cdot (1~8)(\flip{4}~\flip{9})
  \cdot (1~12)(\flip{5}~\flip{10}) \cdot (6~12)(\flip{3}~\flip{9})
  \cdot (0~6)(\flip{4}~\flip{5}) \\
=~ (\flip{1}~\flip{7}~12~6~8)(3~4~\flip{10}~5~\flip{9}). 
\end{multline*}
Here the underlines denote flips;
thus the counter flipped by the move sequence are
$1$, $7$, $9$, and $10$.
Ignoring all the flips
is tantamount to removing all the underlines from the calculation,
which recovers the unsigned permutation of Example~\ref{basic-example}.
\end{example}

Much of the terminology of the previous section
(such as ``closed'', ``degenerate'', etc.)\
carries over to the signed game.
The group of signed permutations of $\PP \sm \{0\}$
induced by closed move sequences is called the
\defterm{signed $\Prp$-game group}, denoted $\GS$;
and the set of signed permutations induced by
the move sequences with $p_0 = 0$
is called $2M_{13}$.

%%%%%%%%%%%%%%%%%%%%%%%%%%%%%%%%%%%%%%%%%%%%%%%%%%%%%%%%%%
%%%%%%%%%%%%%%%%%%%%%%%%%%%%%%%%%%%%%%%%%%%%%%%%%%%%%%%%%%%%%%%%%%%%%%%%
%%%%%%%%%%%%%%%%%%%%%%%%%%%%%%%%%%%%%%%%%%%%%%%%%%%%%%%%%%%%%%%%%%%%%%%%

\section{The signed game, the Golay code, and the Mathieu group}
\label{signed-section}

In this section, we prove the main results
that the basic $\Prp$-game group~$\GM$
is isomorphic to the Mathieu group~$M_{12}$,
and that the signed game group~$\GS$
is the nontrivial double cover~$2M_{12}$.

Let $\Ff_3 = \{0,1,-1\}$ be the field of order~$3$,
and let $X$ be a $13$-dimensional vector space over~$\Ff_3$
with basis $\{x_p \st p \in \PP\}$.
We will write elements of $X$ in the form $v = \sum_p v_p x_p$,
where $v_p \in \Ff_3$.
Define a scalar product on $X$ by
\begin{equation} \label{scalar-product}
v \cdot w = \sum_{p \in \PP} v_p w_p.
\end{equation}
The {\it support\/} of the vector $v$ is
$$
\Supp(v) = \{p \in \PP \st v_p \neq 0\}
$$
and its {\it weight\/} is
$$
\wt(v) = |\Supp(v)|.
$$
We will refer to vector subspaces of $X$ as \defterm{codes},
and to their elements as \defterm{codewords}.
The \defterm{minimal weight} of a code~$X'$ is
$$
\minwt(X') = \min\{\wt(c) \st c \in X', \ c \neq 0\}.
$$

Let $\C \subset X$ be the linear span of the~$13$ vectors
$$
h_{\ell} = \sum_{p \in \ell} x_p
$$
where $\ell$ ranges over $\LL$, and define
$$
\C' = \{c \in \C \ \ \Big\vert \ \ \sum_p c_p = 0\},
$$
a codimension-$1$ subcode of $\C$.
(Note that $\C' \neq \C$ because $h_\ell \notin \C'$.)  
We will show that for each $p \in \PP$,
there is a copy $\C_p$ of the ternary Golay code~\cite[p.~85]{Conway:1999}
occurring naturally as a subcode of $\C$.
First, we set forth some useful properties of $\C$ and $\C'$.

\begin{prop} \label{omnibus-prop} Let $c \in \C$.  Then:
\begin{enumerate}
\item \label{square-eqn}
$\displaystyle \sum_{p \in \PP} c_p^2
= \Biggl(\sum_{p \in \PP}c_p\Biggr)^{\!\!2}$.
\item \label{weight-zero}
$\wt(c) \equiv 0$ or $1 \pmod{3}$.
\item \label{czero-criterion}
$c \in \C'$ iff $\wt(c) \equiv 0 \pmod{3}$.
\item \label{l-identity}
For each $\ell \in \LL$,
$$
\sum_{p \in \PP} c_p = \sum_{p \in \ell} c_p.
$$
\item \label{orthogonal}
$\C' = \C^\bot$, the orthogonal complement of~$\C$
with respect to the scalar product~(\ref{scalar-product}).
\item \label{dimensions}
$\dim\,\C = 7$ and $\dim\,\C' = 6$. 
\item \label{min-weight}
$\minwt(\C)=4$ and $\minwt(\C')=6$.
Moreover, the weight-$4$ codewords in~$\C$
are precisely $\{\pm h_\ell \st \ell \in \LL\}$.
\end{enumerate}
\end{prop}

\begin{proof}
\begin{enumerate}
\item It suffices to show that
\begin{equation} \label{square}
\sum_{p \in \PP} c_p d_p =
\Biggl(\sum_{p \in \PP} c_p \Biggr) \Biggl(\sum_{p \in \PP} d_p \Biggr)
\end{equation}
for all $c,d \in \C$.  Since this identity is bilinear in~$c$ and~$d$,
we need only consider the case $c = h_\ell$, $d = h_m$,
when both sides evaluate to $1$ (whether or not $\ell=m$).

\item
Since the square of each nonzero element of $\Ff_3$ is $1$, we have
$$
\wt(c) \equiv \sum_{p \in \PP} c_p^2 \pmod{3}.
$$
By part~(\ref{square-eqn}),
the right-hand side is a square in $\Ff_3$, hence either $0$ or~$1$.

\item
This follows from the definition of $\C'$,
together with the previous two parts.

\item
It suffices to verify the desired identity for the generators $h_\ell$.
Indeed, let $c = h_m$; then both sides of the identity
are equal to~$1$ whether $\ell$ and~$m$ are the same or different.

\item
If $c \in \C'$, then the right-hand side of (\ref{square})
is zero for every $d \in \C$; it follows that $\C \subset (\C')^\bot$.
To prove the reverse inclusion, let $w \in (\C')^\bot$.
If $\Supp(w)$ intersects some line $\ell$ in more than two points,
then we can reduce $\wt(w)$ by adding $h_\ell$ or $-h_\ell$ to $w$.
Repeating this process, we eventually obtain a codeword
$w' \in C_0^\bot$ which is congruent to $w$ modulo~$\C$
(since $h_\ell \in \C \subset (\C')^\bot$)
and such that $\Supp(w')$ intersects no line in more than two points.
By Proposition~\ref{P3-uniqueness}, $\Prp$ contains no hyperovals,
so $\wt(w') \leq 4$.

Suppose that $\wt(w') \neq 0$.
Then there is a line $\ell$ disjoint from $\Supp(w')$
and another line $m$ intersecting $\Supp(w')$ in exactly one point.
The vector $h_\ell-h_m$ belongs to $C_0$ but is not orthogonal to~$w$,
which is impossible since $\C \subset (\C')^\bot$.
Hence $\wt(w')=0$, $w' = 0$, and $w \in~\C$.

\item
By part (\ref{orthogonal}),
$\dim \C' + \dim (\C')^\bot = 13 = 2\dim \C - 1$.
Hence $\dim \C' = 6$ and $\dim\,\C = 7$.

\item
Clearly $\wt(h_\ell)=4$ for every line $\ell$.
Let $c \in \C$ be a codeword of minimal nonzero weight.
If $\Supp(c)$ meets no line of $\Prp$ in more than two points,
then $c = 0$ by the argument of (\ref{orthogonal}).
In particular $\wt(c) \neq 1$.
By part (\ref{weight-zero}), $\wt(c) \notin \{2,5\}$.
If $\wt(c) \in \{3,4\}$, then $|\Supp(c)  \cap \ell| \geq 3$
for some line $\ell$.
But then the weight of $c$ can be reduced
by adding or subtracting $h_\ell$.
Since $c$ is of minimal weight, this is a contradiction
unless $c = \pm h_\ell$.  Hence $\minwt(\C)=4$.
By part~(\ref{czero-criterion}), we have $\minwt(\C')  \geq 6$.
In fact, $\minwt(\C')=6$ because $\wt(h_\ell-h_m)=6$ for $\ell \neq m$.
\end{enumerate}
\end{proof}

Note that $|\C|=3^7=2187$, which is small enough that
all the assertions of Proposition~\ref{omnibus-prop} could also
be checked by an easily feasible but unenlightening computation.

For each $p \in \PP$, define a subcode
$$
\{c \in \C \st c_p = -\sum_{q \in \PP} c_q\},
$$
and let $\G_p$ be the restriction of $\C_p$
to the coordinates $\PP \sm \{p\}$
(that is, the image of $\C_p$ modulo the subspace spanned by $x_p$).

\begin{prop}
$\G_p$ is isomorphic to the ternary Golay code ${\mathscr C}_{12}$
for every $p \in \PP$.
\end{prop}

\begin{proof}
$\C_p \subsetneq \C$ because $h_\ell \notin \C_p$ for all $\ell$.
The kernel of the restriction map $\phi: \C_p \to \G_p$
can contain only vectors of weight $\leq 1$,
but $\minwt(\C_p)=4$, so $\ker\phi=0$.
Thus $\phi$ is a bijection and $\dim \G_p = \dim \C_p = 6$.

For all $c \in \C_p$,
\begin{eqnarray*}
\wt(\phi(c)) &\equiv& \sum_{q\neq p} c_q^2 \quad \pmod{3} \\
&=& -c_p^2 + \sum_{q \in \PP} c_q^2 \ = \  -c_p^2
  + \Biggl( \sum_{q \in \PP} c_q \Biggr)^{\!\!2} \ \equiv \ 0 \pmod{3};
\end{eqnarray*}
and for all $c,d \in \C_p$,
\begin{eqnarray*}
\phi(c) \cdot \phi(d) &=& \sum_{q \neq p} c_q d_q
  \ = \ -c_p d_p + \sum_{q \in \PP} c_q d_q \\
&=& -c_p d_p + \Biggl( \sum_{q \in \PP} c_q \Biggr)
  \Biggl( \sum_{q \in \PP} d_q \Biggr) \ = \ 0
\end{eqnarray*}
by (\ref{square}).
Hence $\G_p \subseteq \G_p^\perp$, whence $G_p$ is self-dual
  since it has dimension $6=12/2$.
Moreover, $\minwt(\G_p) \geq \minwt(\C) = 4$
  (which implies that $\minwt(\G_p) \geq 6$
  because $\G_p \subseteq \G_p^\perp$).
Therefore $\G_p \isom {\mathscr C}_{12}$~\cite[p.~436]{Conway:1999}.
\end{proof}

Suppose $\ell=\{p,q,r,s\}$.
Let the move $\pmv{p}{q}$ of the signed $\Prp$-game act linearly on~$X$\/
by $\pmv{p}{q}\cdot w = w'$, where
\begin{equation} \label{trdef}
\begin{aligned}
w'_p &= w_q,      \qquad & w'_r &= -w_s, \qquad &
w'_t &= w\0_t ~\text{ for } t \notin \ell, \\
w'_q &= -w_p-w_q, \qquad & w'_s &= -w_r. &  
\end{aligned}
\end{equation}

\begin{prop} \label{golay-iso-prop}
For all $p,q \in \PP$, $\pmv{p}{q}\cdot\C_p\: =\: \C_q$.
\end{prop}

\begin{proof}
Let $\ell=\{p,q,r,s\}$ as above.
Since the linear transformation (\ref{trdef}) is invertible,
it suffices to prove the inclusion $\pmv{p}{q}\cdot\C_p \subset \C_q$.
Let $c \in \C_p$ and $d = \pmv{p}{q}\cdot c$.
By part (\ref{l-identity}) of Proposition~\ref{omnibus-prop},
\begin{eqnarray*}
c_p \ = \  -\sum_{q \in \PP} c_q &=& -\sum_{q \in \ell} c_q \\
&=& -c_p-c_q-c_r-c_s
\end{eqnarray*}
which implies that $c_p = c_q+c_r+c_s$,
since we are working over $\Ff_3$.
Hence
\begin{eqnarray*}
c-d &=& \sum_{p \in \ell} (c_p-d_p)x_p \\
&=& (c\0_p-c\0_q) (x\0_p+x\0_q) + (c\0_r+c\0_s)(x\0_r+x\0_s) \\
&=& (c\0_p-c\0_q) h\0_\ell.
\end{eqnarray*}
So $c-d \in \C$ and $d \in \C$.  Moreover,
\begin{eqnarray*}
\sum_{p \in \PP} d_p \ = \  \sum_{p \in \ell} d_p &=& d_p+d_q+d_r+d_s \\
&=& -c\0_p-c\0_s-c\0_r \\
&=& c\0_p+c\0_q \ = \  -d_q.
\end{eqnarray*}
Therefore $d \in \C_q$.
\end{proof}

A move sequence $[p_0,\dots,p_n]$ acts on~$X$\/
by the composition of the linear transformations~(\ref{trdef})
associated with the moves $\pmv{p_i}{p_{i+1}}$.
It follows from Proposition~\ref{golay-iso-prop}
that the linear transformation associated with $[p_0,\dots,p_n]$
restricts to an isomorphism of $\C_0$ with~$\C_{p_n}$.
In particular, if $p_0 = p_n = 0$,
then $\sigma$ induces an automorphism of the code~$\C_0$.
Accordingly, $\GS$ is naturally isomorphic to a subgroup of $\Aut(\G_0)$.

\begin{prop} \label{sub-mathieu-prop}
$\GM$ is isomorphic to a subgroup of $M_{12}$.
\end{prop}

\begin{proof}
The center $Z$ of $\Aut(\G_0)$ has order two
(it contains the identity map and its negative), 
and $\Aut(\G_0)/Z \isom M_{12}$ (see~\cite[p.~85]{Conway:1999}).
On the other hand, the permutation~$-\Id$
corresponding to the closed path
\begin{equation} \label{minus-id-path}
[0, 10, 7, 0, 4, 1, 2, 4, 3, 5, 6, 3, 0]
\end{equation}
flips each of the~$12$~counters without changing its location.
Clearly $-\Id$ is central and has order~$2$.
Hence $\GS/\{\Id,-\Id\}$ is isomorphic to a subgroup of
$\Aut(\G_0)/Z \isom M_{12}$.
On the other hand, $\GS/\{\Id,-\Id\}=\GM$,
because taking the quotient by $-\Id$ is equivalent to ignoring flips.
\end{proof}

We see now that a permutation $\sigma$ in $M_{12}$ (resp.\ $M_{13}$)
has two \defterm{lifts} $\sigma_1,\sigma_2$
in $2M_{12}$ (resp.\ $2M_{13}$),
both of which are equivalent to $\sigma$ as unsigned permutations
and such that $\sigma_1^{-1} \circ \sigma_2 = -\Id$.

To establish the reverse inclusion, we used a computer program (in~C)
to generate a list of all permutations arising from closed move sequences.\footnote{
%(The source code appears in~\cite{Martin:1996},
%% and may be found online at
%% http://www.math.harvard.edu/$\sim$elkies/mathieu.h}, {\tt 
%% http://www.math.harvard.edu/$\sim$elkies/m13.c}, and {\tt 
%% http://www.math.harvard.edu/$\sim$elkies/2m13.c}.)
%and is now online at
%\verb+<http://www.math.harvard.edu/$\sim$elkies/M13>+.)
  The source code appears in~\cite{Martin:1996},
  and is now online at
  {\tt http://www.math.harvard.edu/$\sim$elkies/M13}~.}  %% end footnote
Rather than reproduce the entire list here,
we use the presentation of $M_{12}$
as a subgroup of the symmetric group~$\Sym_{12}$
with generators given in cycle notation by
\begin{equation} \label{generators}
\begin{aligned}
\alpha &= (1~6~4~2~11~3~8~9~10~7~5), \\
\gamma &= (1~12)(2~9)(3~4)(5~6)(7~8)(10~11), \\
\delta &= (4~5)(2~11)(3~7)(8~9).
\end{aligned}
\end{equation}
(see~\cite[p.~273]{Conway:1999}; we have changed the labelling of the points
to conform with~\eqref{dual-labelling}).
Indeed, the move sequences
$$
[0,11,7,9,8,3,0], \qquad
[0,12,1,9,0,3,8,4,0], \qquad
[0,1,7,0,3,6,0,1,7,0]
$$
induce the permutations $\alpha$, $\gamma$ and $\delta$ respectively.
Combining this with Proposition~\ref{sub-mathieu-prop}
and the known identification of $2M_{12}$ with
$\Aut(\G_0)$~\cite{Conway:1985}, we have proved:

\begin{thm} \label{pthree-mathieu-thm} \quad
\begin{enumerate}
\item The basic $\Prp$-game group $\GM$
is isomorphic to the Mathieu group $M_{12}$,
acting sharply quintuply transitively on $\PP \sm \{0\}$.
\item The signed game group $\GS$
is isomorphic with $2M_{12}$,
with $Z$\/ the $2$-element normal subgroup
and $\GS/Z \isom M_{12}$.
\end{enumerate}
\end{thm}

%%%%%%%%%%%%%%%%%%%%%%%%%%%%%%%%%%%%%%%%%%%%%%%%%%%%%%%%%%%%%%%%%%%%%%%%
%%%%%%%%%%%%%%%%%%%%%%%%%%%%%%%%%%%%%%%%%%%%%%%%%%%%%%%%%%%%%%%%%%%%%%%%
%%%%%%%%%%%%%%%%%%%%%%%%%%%%%%%%%%%%%%%%%%%%%%%%%%%%%%%%%%%%%%%%%%%%%%%%

\section{The dualized game} \label{dualized-section}

We can extend the $\Prp$-game in another way
by placing a second set of counters on the lines of~$\Prp$.
This version of the game provides a second proof
that the game group is $M_{12}$,
realized as the group of automorphisms of a $12 \x 12$
\defterm{Hadamard matrix}
(that is, an orthogonal matrix all of whose entries are $\pm 1$).
In addition, interchanging the roles of points and lines
gives a concrete interpretation of the outer automorphisms of $M_{12}$.

We began the basic $\Prp$-game by placing
$12$ numbered counters on the points $\PP \sm \{0\}$.
In the dualized game, we place in addition
$12$ numbered ``line-counters'' on the lines $\LL \sm \{\pline{0}\}$.
The move sequences of the dualized game are defined
similarly to those of the basic game,
with the proviso that the point-hole must always lie on the line-hole.
Specifically, suppose that the point-hole and line-hole
are located at $p$ and~$\ell$ respectively,
with $\starpt{p} = \{\ell, m, n, k\}$ and $\ell = \{p,q,r,s\}$.
The point-move $\pmv{p}{q}$ is defined as in the basic game;
dually, the line-move $\lmv{\ell}{m}$
consists of moving the line-counter on~$m$ to the hole at~$\ell$
and interchanging the line-counters on~$n$ and~$k$.
Thus a move sequence has the general form
\begin{equation} \label{dualdef}
\left([p_0,\dots,p_n],\:[\ell_0,\dots,\ell_n]\right) \ = \
\lmv{\ell_{n-1}}{\ell_n} \circ \pmv{p_{n-1}}{p_n} \circ\dots\circ
\lmv{\ell_0}{\ell_1} \circ \pmv{p_0}{p_1}
\end{equation}
subject to the conditions
$p_i,\,p_{i+1} \in \ell_i$ and
$\ell_i,\,\ell_{i+1} \in \starpt{p_{i+1}}$ for all $i$.
Each move sequence induces a pair of permutations
$\sigma = (\sigma_\PP,\sigma_\LL)$,
where $\sigma_{\PP}$ acts on the point-counters
and $\sigma_{\LL}$ acts on the line-counters.

It is easy to verify that for every move sequence
$\left([p_0,\dots,p_n],\:[\ell_0,\dots,\ell_n]\right)$
of the dualized game, the point-path $[p_0,\dots,p_n]$ is nondegenerate
if and only if the line-path $[\ell_0,\dots,\ell_n]$ is.  
As before, every move sequence
is equivalent to one in which both paths are nondegenerate.

A move sequence of the dualized game is called \defterm{closed}\/
if it returns both the point-hole and the line-hole
to their initial locations.
The group of permutations induced by closed moves is called the
\defterm{dualized $\Prp$-game group}, written~$\GD$.
In fact, we shall show that $\GD \isom \GM$,
and indeed that the point-permutation of an element of~$\GD$
determines the line-permutation uniquely and vice versa.

\begin{example} \label{dualized-example}
As in Examples~\ref{basic-example} and~\ref{signed-example}),
consider the path $[0,6,12,1,8,0]$.
For this to be the point-path
of a closed move sequence in the dualized game,
the corresponding line-path can only be
$$
\left[\twopt{0}{\ 6},\ \ \twopt{6}{\ 12},\ \ \twopt{12}{\ 1},\ \ 
\twopt{1}{\ 8},\ \ \twopt{8}{\ 0}\right] \ = \ 
\left[ \pline{0}, \pline{1}, \pline{12}, \pline{6}, \pline{4} \right].
$$
\end{example}

The moves of the dualized game
may be interpreted as automorphisms of a $12 \x 12$ Hadamard matrix~$H$.
An \defterm{automorphism} of~$H$\/
may be defined as a pair $(\sigma,\tau)$ of signed permutation matrices
such that $\sigma H \tau = H$.
The group $\Aut(H)$ of all automorphisms
is isomorphic to $2M_{12}$~\cite[p.~32]{Conway:1985}.
In what follows, we construct an isomorphism of~$\GD$ with~$\Aut(H)$.

Define a modified incidence matrix $E = (e_{ij})$ for $\Prp$,
with rows $i$ indexed by $\PP$ and columns $j$ indexed by $\LL$, by
$$
e_{ij} = \begin{cases}
-1  & \quad i \in j, \\
+1  & \quad i \notin j.
\end{cases}
$$
Labelling the points and lines of $\Prp$ self-dually,
as in~(\ref{dual-labelling}), makes $E$\/ into a symmetric matrix.
Each row of $E$ contains four $-1$'s and nine $+1$'s,
and each pair of distinct rows agree in exactly seven columns,
so the scalar product $E_r \cdot E_s$ of two rows of~$E$\/ is
\begin{equation} \label{scalar-product-formula}
E_r \cdot E_s \ = \  \sum_j e_{rj} e_{sj} \ = \  \begin{cases}
13 & \quad r=s, \\
1 & \quad r \neq s.
\end{cases}
\end{equation}
Next, for all pairs $p,\ell$ with $p \in \ell$,
define a $12 \x 12$ matrix $\Hadmat{p}{\ell}$,
with rows indexed by $\PP \sm \{p\}$
and columns indexed by $\LL \sm \{\ell\}$, by
\begin{equation} \label{define-H_pl}
(\Hadmat{p}{\ell})_{ij} \ = \  \begin{cases}
-e_{ij} & \quad i \in \ell ~\text{and}~ j \in \starpt{p} \\
e_{ij} & \quad \otherwise.
\end{cases}
\end{equation}

\begin{prop}
Let $p \in \PP$ and $\ell \in \starpt{p}$.
Then $H=\Hadmat{p}{\ell}$ is a Hadamard matrix.
\end{prop}

\begin{proof} 
$H$ can be made symmetric by choosing a self-dual labelling
in which $p$ and $\ell$ have the same label.
Thus, to prove the proposition, it is enough to show that 
for each pair of distinct points $r,s \in \PP\sm\{p\}$,
the scalar product $H_r \cdot H_s$ of the corresponding rows of~$H$ is zero.
%, that is,
%	\begin{equation} \label{scalar-prod-H}
%	H_r \cdot H_s \ = \ \sum_{j \in \LL \sm \{\ell\}} e_{rj} e_{sj} \ = \ 0.
%	\end{equation}

If $r \notin \ell$ and $s \notin \ell$, then $e_{r\ell}=e_{s\ell}=1$, so
$$
H_r \cdot H_s \ = \ \sum_{j \in \LL \sm \{\ell\}}^{\phantom{0}} e_{rj} e_{sj} \ = \
-e_{rl} e_{sl}\;+\;\sum_{j \in \LL} e_{rj} e_{sj}
\ = \ -1 +  E_r \cdot E_s \ = \ 0.
$$
If $r \in \ell$ and $s \in \ell$, then $e_{r\ell}=e_{s\ell}=-1$, so
$$
H_r \cdot H_s \ = \ \sum_{j \in \starpt{p} \sm \{\ell\}} (-e_{rj}) (-e_{sj})
\ + \sum_{j \in \LL \sm \starpt{p}} e_{rj} e_{sj}
\ = \ \sum_{j \in \LL \sm \{\ell\}}^{\phantom{0}} e_{rj} e_{sj} \ = \ 0
$$
by the previous case.
Finally, suppose that $r \in \ell$ and $s \notin \ell$.
Then $e_{r\ell} = -1$ and $e_{s\ell} = 1$, so
\begin{equation*}
H_r \cdot H_s \ = \ 
-\sum_{j \in \starpt{p} \sm \{\ell\}} e_{rj} e_{sj} ~+~
\sum_{j \in \LL \sm \starpt{p}} e_{rj} e_{sj}.
\end{equation*}
Moreover,
\begin{equation*}
1 \ = \ \sum_{j \in \LL} e_{rj} e_{sj}
\ = \  -1 + \sum_{j \in \LL \sm \starpt{p}}
e_{rj} e_{sj} ~+ \sum_{j \in \starpt{p} \sm \{\ell\}} e_{rj} e_{sj}.
\end{equation*}
by (\ref{scalar-product-formula}).
Combining these two observations, we obtain
\begin{equation}\label{hadform}
H_r \cdot H_s \ = \ 2 \left(1-\sum_{j \in \starpt{p}\sm\{\ell\}} e_{rj} e_{sj} 
\right).
\end{equation}
Since $\twopt{p}{r} = \ell$,
the three $e_{rj}$'s in the right-hand side of (\ref{hadform})
all equal $+1$.
On the other hand, $\twopt{p}{s} \neq \ell$,
so $\twopt{p}{s}$ is one of the other lines in $\starpt{p}$.
Thus one of the three $e_{sj}$'s is $-1$ and the other two are $+1$.
Therefore the expression in (\ref{hadform}) vanishes.
\end{proof}   

We now associate a signed permutation matrix
with each move sequence of the dualized game.  
For $p \in \PP$ and $\ell = \{p,q,r,s\} \in \starpt{p}$,
let $B = (b_{ij})$ be a $12 \x 12$ matrix,
with rows indexed by $\PP \sm \{p\}$
and columns indexed by $\LL \sm \{\ell\}$.
The point-move $\pmv{p}{q}$ acts on $B$,
producing a matrix $\pmv{p}{q}\cdot B$\/
with rows indexed by $\PP \sm \{q\}$
and columns indexed by $\LL \sm \{\ell\}$,
whose $(i,j)$ entry is
\begin{equation} \label{act-hadamard}
\left(\pmv{p}{q} \cdot B\right)_{ij} \ = \  \left\{ \begin{array}{rl}
 b_{qj} ~&~ i=p \\
-b_{rj} ~&~ i=s \\
-b_{sj} ~&~ i=r \\
 b_{ij} ~&~ \otherwise.
\end{array} \right.
\end{equation}
The line-move $\lmv{\ell}{m}$ acts on the columns of~$B$\/
in a similar way, producing a matrix
with rows indexed by $\PP \sm \{p\}$
and columns indexed by $\LL \sm \{m\}$.
More generally, we may associate a signed permutation matrix
with each move sequence of the dual game
by composing those corresponding to its constituent moves.
Note that the actions of point- and line-moves commute.

\begin{prop} \label{Hadamard-transform}
Let $\sigma = (\sigma_\PP,\sigma_\LL)$
be a move sequence of the dualized game, with the point-hole initially at $p\in\PP$
and the line-hole initially at $\ell\in\LL(p)$.
Let $\Hadmat{p}{\ell}$ be the Hadamard matrix defined in \eqref{define-H_pl}.
Then $\sigma(\Hadmat{p}{\ell}) = \Hadmat{\sigma(p)}{\sigma(\ell)}$.
\end{prop}

%\begin{prop} \label{Hadamard-transform}
%Let $p \in \ell$ and $q \in m$.
%Define Hadamard matrices $\Hadmat{p}{\ell}$ and
%$\Hadmat{q}{m}$ as in (\ref{define-H_pl}).
%Let $\sigma = (\sigma_\PP,\sigma_\LL)$
%be a move sequence of the dualized game, with
%$\sigma_\PP(p) = q$ and $\sigma_\LL(\ell) = m$.  
%Then $\sigma(\Hadmat{p}{\ell}) = \Hadmat{q}{m}$.
%\end{prop}

\begin{proof} It is sufficient to consider the case that
$\sigma$ is a single point-move.  The proof for line-moves is identical,
and the general case will then follow by composition.
Suppose therefore that $\ell = \{p,q,r,s\}$ and $\sigma = \pmv{p}{q}$.
The definition (\ref{define-H_pl}) may be rewritten as
$$
(\Hadmat{p}{\ell})_{ij} \ = \  \left\{ \begin{array}{rlll}
-e_{ij} ~&~ \quad i \in \{q,r,s\} & \text{and} & j \in \starpt{p} \sm \{\ell\} \\
+e_{ij} ~&~ \quad i \in \PP \sm \ell & \text{and} & j \in \LL \sm \starpt{p},
\end{array} \right.
$$
so that
$$
(\sigma(\Hadmat{p}{\ell}))_{ij} \ = \  \begin{cases}
(\Hadmat{p}{\ell})_{qj} ~&~ \quad i=p \\
-(\Hadmat{p}{\ell})_{sj} ~&~ \quad i=r \\
-(\Hadmat{p}{\ell})_{rj} ~&~ \quad i=s \\
(\Hadmat{p}{\ell})_{ij} ~&~ \quad \otherwise,
\end{cases}
$$
and
$$
(\Hadmat{q}{\ell})_{ij} \ = \  \begin{cases}
-e_{ij} & \quad i \in \ell ~\text{and}~ j \in \starpt{q} \\
+e_{ij} & \quad \otherwise.
\end{cases}
$$

We will show that
$(\sigma(\Hadmat{p}{\ell}))_{ij} = (\Hadmat{q}{\ell})_{ij}$
for all $i,j$.  
First, if $i \notin \ell$, then
$$
(\sigma(\Hadmat{p}{\ell}))_{ij} \ = \ 
(\Hadmat{p}{\ell})_{ij} \ = \  e_{ij} \ = \  
(\Hadmat{q}{\ell})_{ij}.
$$
Second, suppose that $i=p$.  In this case
\begin{eqnarray*}
(\sigma(\Hadmat{p}{\ell}))_{pj} \ = \ 
(\Hadmat{p}{\ell})_{qj} &=&
\begin{cases}
-e_{qj} & \quad j \in \LL \sm \starpt{p} \\
+e_{qj} & \quad \otherwise
\end{cases} \\
&=& \begin{cases}
-1 & \quad j \in \starpt{p} \cup \starpt{q} \sm \{\ell\} \\
+1 & \quad \otherwise
\end{cases} \\
&=& \begin{cases}
-e_{pj} & \quad j \in \starpt{q} \sm \{\ell\} \\
+e_{pj} & \quad \otherwise
\end{cases} \\
&=& (\Hadmat{q}{\ell})_{pj}.
\end{eqnarray*}
Finally, suppose that $i=r$ (the case $i=s$ is analogous).  Then
\begin{eqnarray*}
(\sigma(\Hadmat{p}{\ell}))_{rj} \ = \  -(\Hadmat{p}{\ell})_{sj} &=&
\begin{cases}
+e_{sj} & \quad j \in \starpt{p} \sm \{\ell\} \\
-e_{sj} & \quad \otherwise
\end{cases} \\
&=& \begin{cases}
+1 & \quad j \in \starpt{p} \cup \starpt{s} \sm \{\ell\} \\
-1 & \quad j \in \starpt{q} \cup \starpt{r} \sm \{\ell\} \\
\end{cases} \\
&=& \begin{cases}
-e_{rj} & \quad j \in \starpt{q} \sm \{\ell\} \\
+e_{rj} & \quad \otherwise
\end{cases} \\
&=& (\Hadmat{q}{\ell})_{rj}.
\end{eqnarray*}
\end{proof}   

\begin{cor}
$\GD\isom\GM$ {\rm (}$\isom M_{12}${\rm )}.
\end{cor}

\begin{proof}
Proposition~\ref{Hadamard-transform} implies that for each closed move $\sigma$ of
the dualized game, the pair of (unsigned) permutations $(\sigma_\PP,\sigma_\LL)$ is
an automorphism of the Hadamard matrix $H=\Hadmat{p}{\ell}$.
That is, we have an injective group homomorphism from $\GD$ to $\Aut(H)/\{\pm\Id\}\isom M_{12}$.
On the other hand, the permutations \eqref{generators} generate a subgroup of $\GD$ that
is isomorphic to $M_{12}$.
\end{proof}

Denote by $\Aut(M_{12})$ the group of automorphisms of $M_{12}$,
and by $\Inn(M_{12})$ the normal subgroup of inner automorphisms
(that is, automorphisms given by conjugation).
Then $\Inn(M_{12}) \isom M_{12}$
since $M_{12}$ is simple and nonabelian,
and the quotient $\Aut(M_{12}) / \Inn(M_{12})$
has order two~\cite[p.~31]{Conway:1985}.
The dualized game allows us to describe an outer automorphism
(and thus the full automorphism group) of~$M_{12}$ explicitly.
Consider the map
\begin{equation} \label{define-theta}
\begin{array}{cccc}
\theta: & \GD & \to & \GD \\
& (\sigma_P,\sigma_L) & \mapsto & (\sigma_L,\sigma_P).
\end{array}
\end{equation}

\begin{prop}
The map $\theta$ is an outer automorphism of $\GD \isom M_{12}$.
\end{prop}

\begin{proof}
The map $\theta$ respects concatenation of paths,
so it is a group homomorphism $\GD \to \GD$.
It is clearly surjective, hence an automorphism.  It remains only
to show that $\theta$ is not conjugation by any element
of $\GD$.

Consider the point-paths $\pi_1=[0,1,4,0]$, $\pi_2=[0,2,10,0]$ 
and $\pi_3=[0,3,12,0]$, whose induced permutations are respectively
  $$\alpha_1 = (1\ 4)(2\ 3)(5\ 6)(8\ 9),  \quad
    \alpha_2 = (1\ 3)(2\ 10)(6\ 8)(9\ 11),\quad
    \alpha_3 = (1\ 2)(3\ 12)(6\ 9)(7\ 8).$$
By the labelling \eqref{dual-labelling}, for each $i$, the line-path
corresponding to $\pi_i$ is its reverse, so $\theta(\alpha_i)=\alpha_i^{-1}=\alpha_i$
(because each $\alpha_i$ is an involution).
Therefore, if $\theta$ is conjugation by some permutation $\sigma$, then
$\sigma$ must commute with each $\alpha_i$.
A Maple computation (which is not hard to check by hand) reveals that the
intersection of the \hbox{$\Sym_{12}$-centralizers} of the $\alpha_i$'s
contains exactly one non-identity element, namely
  $$\sigma = (1\ 6)(2\ 9)(3\ 8)(4\ 5)(7\ 12)(10\ 11).$$
(In fact, this permutation commutes with every fixed point
of the automorphism~$\theta$.)
Now, consider the point-path $[0,1,5,0]$, whose associated line-path is $[0,6,1,0]$.
The induced point- and line-permutations are respectively
$(1\ 5)(2\ 3)(4\ 6)(10\ 12)$ and $(1\ 6)(2\ 3)(4\ 5)(7\ 11)$.
Then $\theta$ interchanges these two permutations; however, they are not
conjugates under $\sigma$.  It follows that $\theta$ is not conjugation by
any element of $\Sym_{12}$, so {\em a~fortiori} not by any element
of $M_{12}$; that is, $\theta$ is an outer automorphism.
\end{proof}

\begin{cor} $\langle \GD,\theta \rangle = \Aut(\GD) = \Aut(M_{12})$.
\end{cor}

%%%%%%%%%%%%%%%%%%%%%%%%%%%%%%%%%%%%%%%%%%%%%%%%%%%%%%%%%%%%%%%%%%%%%%%%
%%%%%%%%%%%%%%%%%%%%%%%%%%%%%%%%%%%%%%%%%%%%%%%%%%%%%%%%%%%%%%%%%%%%%%%%
%%%%%%%%%%%%%%%%%%%%%%%%%%%%%%%%%%%%%%%%%%%%%%%%%%%%%%%%%%%%%%%%%%%%%%%%

\section{$M_{13}$ and sextuple transitivity} \label{sextuple-section}
\subsection{Multiply transitive groups}

We recall some basic terminology.
Let $G$ be a (finite) group acting on a (finite) set $X$\/;  
that is, there is a group homomorphism from $G$ to $\Sym_X$,
the group of permutations of $X$.  
The action is called \defterm{faithful}
if this homomorphism is one-to-one and
\defterm{transitive} if the action has a single orbit.
More generally, the action is \defterm{$k$-transitive}
if $|X| \geq k$ and for any two $k$-tuples of distinct elements of~$X$,
say $\ptuple=(p_1,\dots,p_k)$ and $\qtuple=(q_1,\dots,q_k)$,
there exists $g \in G$ such that $g \cdot p_i = q_i$ for all $i$.
If the element $g$ is unique,
then the action is \defterm{sharply $k$-transitive}.
Note that a group $G$\/ with a faithful,
sharply \hbox{$k$-transitive} action on an $r$-set
must have cardinality $r!/(r-k)!$.

Groups with highly transitive actions are quite unusual:
by the classification of finite simple groups,
the only groups with a sharply quintuply transitive action are
$M_{12}$, $\Sym_5$, $\Sym_6$ and $\Alt_7$,
  and there are no sextuply transitive groups
  other than $\Sym_7$ and $\Alt_8$.
In particular---and this does not require
the Classification Theorem---it is not possible to continue
Mathieu's construction of $M_{11}$ and~$M_{12}$
to a transitive subgroup of~$\Sym_{13}$
other than $\Alt_{13}$ and $\Sym_{13}$ itself.
Yet we have obtained the pseudogroup $M_{13}$
by a method similar to this construction.
$M_{12}$, realized as the game group $\GM$,
acts faithfully and sharply \hbox{$5$-transitively}
on the $12$-element set $\PP \sm \{0\}$.
Meanwhile, $M_{13}$ ``acts'' on the $13$-element set~$\PP$,
and $|M_{13}|=13|M_{12}|$,
just what the order of a sharply sextuply transitive group ought to be.

We are thus led to consider the question:
Is the ``action'' of $M_{13}$ on $\PP$ sextuply transitive?
That is, given two sextuples
$\ptuple = (p_1,\dots,p_6)$ and $\qtuple = (q_1,\dots,q_6)$
of points of $\Prp$, does there exist some $\sigma \in M_{13}$
such that $\sigma(p_i)=q_i$ for all $i$?
(Here and from now on, ``sextuple''
means ``ordered sextuple of distinct elements''.
In addition, we wish to include the possibility that $0 \in \ptuple$,
so ``counter'' really means ``counter or hole''.)
Since $M_{13}$ is not a group, there are actually two distinct questions:

\begin{enumerate}
\item
Fix a sextuple $\ptuple$ of counters.
Is it true that for all sextuples $\qtuple$ of points of $\Prp$,
there exists some $\sigma \in M_{13}$
such that $\sigma(\ptuple) = \qtuple$?
If so, we call $\ptuple$ a \defterm{universal donor}.
\item
Fix a sextuple $\qtuple$ of points of $\Prp$.
Is it true that for all sextuples $\ptuple$ of counters,
there exists some $\sigma \in M_{13}$
such that $\sigma(\ptuple)=\qtuple$?
If so, we call $\ptuple$ a \defterm{universal recipient}.
\end{enumerate}

We will examine the questions separately.
In each case, our computational data was invaluable
as a source of educated guesses about sextuple transitivity.
We start by making an elementary observation
which will be quite useful in both cases.

\begin{lemma}\label{pigeon}
Let $\ptuple$ be a sextuple of counters.
Then $\ptuple$ is a universal donor if and only if,
for all $\sigma,\tau \in M_{13}$ with $\sigma \neq \tau$,
we have $\sigma(\ptuple) \neq \tau(\ptuple)$.
Similarly, if $\qtuple$ is a sextuple of points,
it is a universal recipient if and only if
$\sigma \neq \tau$ implies $\sigma^{-1}(\qtuple) \neq \sigma^{-1}(\qtuple)$.
\end{lemma}

\begin{proof}
This follows from the pigeonhole principle,
together with the observation that $|M_{13}|$ equals $13!/7!$,
the number of sextuples of points in $\Prp$.
\end{proof}   

%%%%%%%%%%%%%%%%%%%%%%%%%%%%%%%%%%%%%%%%%%%%%%%%%%%%%%%%%%

\subsection{Sextuple transitivity on counters}

We consider the question of when a sextuple $\ptuple$ of counters
is a universal donor.  Note that the property is invariant
under permuting the order of the $p_i$.
Thus, for ease of notation, we frequently treat $\ptuple$ as a set:
for instance, we write $\ptuple \cap \ell$
rather than $\{p_i \st 1 \leq i \leq 6\} \cap \ell$.

\begin{thm} \label{UD-thm}
A sextuple of counters $\ptuple = (p_1, \ldots, p_6)$ is
a universal donor if and only if $p_i=0$ for some $i$.
\end{thm}

\begin{proof}
Suppose first that $p_i=0$ for some $i$.  Let $\qtuple = (q_1,\dots,q_6)$ be an
arbitrary sextuple of points.  Note that the move $\pmv{0}{q_i}$ takes the hole
from $p_i$ to~$q_i$.  The $q_i$-conjugate of~$\GM$ acts quintuply transitively
on $\PP \sm \{q_i\}$, hence contains a permutation $\sigma$ such that
$$
\big( \sigma \;\circ\; \pmv{0}{q_i} \big)\, (p_j) \ = \  q_j
$$
for all $j \neq i$.  That is,
$\sigma \circ \pmv{0}{q_i}$ is the desired element of~$M_{13}$
taking $\ptuple$ to $\qtuple$.

Now suppose that $0\not\in \ptuple$.  The set $\PP-\ptuple-\{0\}$ has cardinality six.
Since $\Prp$ has no hyperhyperovals (as discussed in Section~\ref{pthree-subsection}), there is some line $\ell$
that meets $\PP-\ptuple-\{0\}$ in at least three points; that is, the set $A = (\ptuple\cup\{0\})
\cap \ell$ has at most one element.  We consider three cases;
in each case, we will exhibit two moves $\sigma,\tau\in M_{13}$ that act equally on the counters
of $\ptuple$; by Lemma~\ref{pigeon}, such a pair will suffice to show that $\ptuple$ is not a universal donor.

%in each case, we will construct a move $\tau \in M_{13} \sm \{\Id\}$
%such that $\tau(\ptuple)=\ptuple$.  Together with Lemma~\ref{pigeon}, this
%will imply that $\ptuple$ is not a universal donor.

\Case{1}{$A = \{0\}$}.
Then $\ell \cap \ptuple = \emptyset$,
so $\pmv{0}{q}$ fixes each counter in $\ptuple$ for any point $q\in\ell$ other than $0$.
Thus we may take $\sigma=\Id$ and $\tau = \pmv{0}{q}$.

\Case{2}{$A = \{p_i\}$ for some $i$}.  Let $q$ be any point on $\ell$ other than $p_i$.
Playing the move $[0,p_i]$ results in a position in which $\ell$ contains no counters
of $\ptuple$; therefore, we may take $\sigma=[0,p_i]$ and $\tau=[0,p_i,q]$.

%\Case{2}{$A = \{p_i\}$ for some $i$}.  Then we may take $\tau = \pmv{0}{p_i}$.
%\Case{3}{$A = \emptyset$}.  Then we may take $\tau = \pmv{0}{r}$ for any $r \in \ell$.

\Case{3}{$A = \emptyset$}.  Let $q,r$ be distinct points on $\ell$.
Similarly to Case~2, we may take $\sigma=[0,q]$ and $\tau=[0,q,r]$.
\end{proof}

%%%%%%%%%%%%%%%%%%%%%%%%%%%%%%%%%%%%%%%%%%%%%%%%%%%%%%%%%%%

\subsection{Sextuple transitivity on points}

We now consider the question
of when a sextuple $\qtuple$ of points is a universal recipient.
As before, we shall make no notational distinction
between the ordered sextuple $\qtuple$ and its underlying set.

\begin{thm} \label{UR-thm}
A sextuple of points $\qtuple = (q_1,\dots,q_6)$
is a universal recipient if and only if it contains some line of $\Prp$.
\end{thm}

\begin{proof} Suppose that $\qtuple$ contains a line $\ell$.
Let $\ptuple$ be a sextuple of counters;
our goal is to find $\sigma \in M_{13}$ taking $\ptuple$ to $\qtuple$.
If $0 \in \ptuple$, then $\ptuple$ is a universal donor
by Theorem~\ref{UD-thm} so we are done.
Suppose now that $0 \notin \ptuple$.
Without loss of generality we may suppose
that $\ell = \{q_1,q_2,q_3,q_4\}$,
and that the line $m = \twopt{q_5}{q_6}$ meets $\ell$ at $q_1$.
Let $x$ be the fourth point on this line.

By quintuple transitivity, the $q_1$-conjugate of~$\GM$
contains a move $\tau$ such that $\tau(p_i)=q_i$ for $2 \leq i \leq 6$.
Consider the move $\upsilon = \tau \circ \pmv{0}{q_1}$;
note that $\upsilon(p_1) \notin \qtuple$.
If $\upsilon(p_1) \neq x$, then the move
$\pmv{q_1}{\upsilon(p_1)} \circ \upsilon$
moves the counter $p_1$ to the point $q_1$
but does not move any other counter in $\ptuple$.
Hence the desired permutation $\sigma \in M_{13}$
taking $\ptuple$ to $\qtuple$ is
$$
\sigma \ = \ \pmv{q_1}{\upsilon(p_1)} \;\circ\; \tau \;\circ\; \pmv{0}{q_1}.
$$

On the other hand, suppose $\upsilon(p_1) = x$.
The move sequence $[q_1,q_2,x,q_3,q_4,x]$ induces the permutation
$$
\rho \ = \  (r_2~s_2)(r_3~s_3)(r_4~s_4)(q_1~x),
$$
where $\twopt{x}{q_i} = \{x,q_i,r_i,s_i\}$ for $i=2,3,4$.
Hence the move
$$
\sigma \ = \  \rho \;\circ\; \upsilon \;\in\; M_{13}
$$
takes $p_i$ to $q_i$ for all $i$ as desired.

For the ``only if'' direction of the theorem, suppose that $\qtuple$ does not contain any line.
We will show that there are two distinct elements of~$M_{13}$
which carry the same ordered sextuple of counters to the points~$q_i$.
It will follow by Lemma~\ref{pigeon} that $\qtuple$ is not a
universal recipient.

If $\ell \cap \qtuple = \emptyset$ for some line $\ell$,
then our task is easy.  Let $p_1,p_2 \in \ell$.
Then $\pmv{0}{p_2}$ and $\pmv{p_1}{p_2} \circ \pmv{0}{p_1}$
are elements of $M_{13}$ carrying the same set of counters to $\qtuple$.
By Lemma~\ref{pigeon}, $\qtuple$ is not a universal recipient.

The more difficult case is when $\qtuple$ meets every line,
but does not contain any line.  Since $\Prp$ has no hyperhyperovals,
%By the pigeonhole principle, two of the~$15$ lines $\twopt{q_i}{q_j}$
%for $i<j$ must be equal; thus
we may assume that $q_1,q_2,q_5$ lie on a common line $\ell$
(the reason for this apparently strange choice will be clear momentarily).
Let $y$ be the fourth point on $\ell$.
Then $y \notin \qtuple$, and for $\qtuple$ to meet every line,
each of the points $q_3,q_4,q_6$
must lie on a different line in $\LL(y) \sm \{\ell\}$.
Thus each of the lines
$\twopt{q_3}{q_6}$, $\twopt{q_4}{q_6}$, $\twopt{q_3}{q_4}$
meets $\ell$ in a point other than $y$.
Without loss of generality we may assume that
$$
q_1 \in \twopt{q_3}{q_6}, \qquad q_2 \in \twopt{q_4}{q_6}, \qquad
q_5 \in \twopt{q_3}{q_4}.
$$
In particular, the points $q_1,q_2,q_3,q_4$ form an oval.
Thus we may adopt the labelling \eqref{oval-labelling}, with
$$
q_5 = \twopt{q_1}{q_2} \cap \twopt{q_3}{q_4} = r_1, \qquad
q_6 = \twopt{q_1}{q_3} \cap \twopt{q_2}{q_4} = r_2, \qquad
y = s_{12}.
$$

%We now construct two elements of $M_{13}$ which bring the same counters
%to $\qtuple = (q_1,\dots,q_6) = (q_1,q_2,q_3,q_4,r_1,r_2)$.
%First
If $s_{12}=0$, then the paths $[s_{12},r_2,s_{23}]$ and $[s_{12},r_3,r_1,s_{14}]$
induce the permutations
  \begin{align*}
  \sigma &= (s_{12}\ s_{23}\ r_2) (r_3\ s_{34}) (r_1\ s_{14}),\\
  \tau &= (s_{12}\ s_{14}\ r_1\ r_3) (r_2\ s_{34}\ s_{23}) (s_{13}\ s_{24})
  \end{align*}
respectively.  Both of these elements of $M_{13}$ fix $q_1,q_2,q_3,q_4$,
and move the counters originally located at $s_{14},s_{23}$
to $r_1=q_5$ and $r_2=q_6$ respectively.
Therefore, by Lemma~\ref{pigeon}, $\qtuple$ is not a universal recipient.
On the other hand, if $s_{12}\neq 0$, then we need only preface the moves
$\sigma,\tau$ given above by moving the hole to $s_{12}$.
That is, the moves $[0,s_{12},r_2,s_{23}]$ and $[0,s_{12},r_3,r_1,s_{14}]$
move the same ordered sextuple of counters to the points $\qtuple$.
\end{proof}

%%%%%%%%%%%%%%%%%%%%%%%%%%%%%%%%%%%%%%%%%%%%%%%%%%%%%%%%%%%%%%%%%%%%%%%%
%%%%%%%%%%%%%%%%%%%%%%%%%%%%%%%%%%%%%%%%%%%%%%%%%%%%%%%%%%%%%%%%%%%%%%%%
%%%%%%%%%%%%%%%%%%%%%%%%%%%%%%%%%%%%%%%%%%%%%%%%%%%%%%%%%%%%%%%%%%%%%%%%

\section{Metric properties} \label{metric-section}
\subsection{The basic game}

Let $G$\/ be a group generated by a finite set $X$.
The \defterm{Cayley graph} of~$G$\/ with respect to~$X$\/
is the graph whose vertices are the elements of~$G$,
with $g,g'$ connected by an edge if $g=xg'$ for some $x \in X$.
We define the Cayley graph $\Cayley$ of $M_{13}$ analogously:
the vertices are the $13!/7!$ positions of the basic game,
and two positions are connected by an edge
if one may be obtained from the other by a single move $\pmv{p}{q}$.

We may use the Cayley graph to define a metric on $M_{13}$,
as follows: $d(\sigma,\tau)$
is the length of the shortest path in $\Cayley$
with endpoints $\sigma$ and $\tau$, that is,
the minimal number of moves needed to go from $\sigma$ to $\tau$.
Note that no two elements of $M_{12}$ are adjacent in $\Cayley$.
Indeed, $d(\sigma,\tau) \geq 3$ for $\sigma \neq \tau \in M_{12}$,
because a two-move path returning the hole to the starting position
must be of the form $\pmv{p}{q} \circ \pmv{q}{p} = \Id$.
Also, a path from $\sigma$ to $\tau$ with length exactly $d(\sigma,\tau)$
must be nondegenerate.  The \defterm{depth} of $\sigma$
is defined as $d(\sigma) = d(\sigma,\Id)$.
We also define
\begin{equation} \label{define-numd}
\begin{aligned}
\numd{M_{12}}{k} ~&=~ \#\left\{\sigma \in M_{12} \st d(\sigma)=k\right\}, \\
\numd{M_{13}}{k} ~&=~ \#\left\{\sigma \in M_{13} \st d(\sigma)=k\right\}.
\end{aligned}
\end{equation}

We can find these numbers from the computer-generated table of move sequences.

\begin{prop} \label{depth-unsigned-prop}
The depth distributions for $M_{12}$ and $M_{13}$
are given by the following table:
$$
\begin{array}{c|rccccccccc}
k & 0 & 1 & 2 & 3 & 4 & 5 & 6 & 7 & 8 & 9 \\ \hline
\numd{M_{12}}{k} & 1 &  0 &   0 &  54 &  540 &  5184 &  25173 &  55044 &  9036 & 8 \\
\numd{M_{13}}{k} & 1 & 12 & 108 & 918 & 7344 & 57852 & 344925 & 733500 & 90852 & 8
\end{array}
$$
\end{prop}

We can explain some of the ``shallower'' numbers
without resorting to computation.
The unique element at depth~$0$ is obviously the identity.
There are no moves in $M_{12}$ at depths~$1$ or~$2$
because there are no nondegenerate closed paths having those lengths.

Let $[0,p,q,0]$ be a nondegenerate closed path of length~$3$.
For nondegeneracy, we must have $p \neq 0$ and $q \notin \twopt{0}{p}$,
so there are $12 \cdot 9 = 108$ such paths.
The permutation induced by each path has cycle-shape $2^4$
(that is, it is a quadruple transposition).
This permutation has order~$2$, so the path $[0,q,p,0]$ is equivalent.
This is the reason that $\numd{M_{12}}{3}=108/2=54$.

Let $[0,p,q,r,0]$ be a nondegenerate closed path of length~$4$.
For nondegeneracy, we must have $p \neq 0$, $q \notin \twopt{0}{p}$,
and $r \notin \twopt{p}{q} \cup \twopt{q}{0}$,
so there are $12 \cdot 9 \cdot 6 = 648$ such paths.
If $\{0,p,r\}$ are collinear,
then the cycle-shape of the induced permutation is $2^4$;
otherwise it is $3^3$.
In the first case, the path $[0,r,q,p,0]$ is equivalent.
There are $216$ paths with $\{0,p,r\}$ collinear,
so $108$ of them are redundant.
Since $648-108=540=\numd{M_{12}}{4}$,
there are no other equivalences among paths of this length.

The computer data may also be used to tabulate
the nondegenerate closed paths of length~$k$
inducing the identity permutation.
There are no such paths of length $k<6$,
and the paths of length $k=6,7,8$
are unique up to automorphisms of $\Prp$.
For $k=6$, all such paths have the form
$$
[0,p,q,0,p,q,0]
$$
where $0,p,q$ are noncollinear.  For $k=7$, all paths have the form
$$
[0,p,r,q,p,r,q,0]
$$
where $0,p,q$ are collinear and $r$ does not lie on their common line.
For $k=8$, all paths have the form
$$
[0,p,q,r,p,q,r,p,0]
$$
where $\{0,p,q,r\}$ is an oval.
In particular, the number of length-$8$ paths inducing the identity
is the number of ordered ovals beginning with $0$,
which is $12 \cdot 9 \cdot 4 = 432$.
Note that by Proposition~\ref{P3-uniqueness}, this is the cardinality
of the stabilizer of a point in $\Aut(\Prp) = \PGL_2(\Ff_3)$.

A striking feature of the depth distribution
is that there are only eight permutations at maximal depth.
These permutations are
\begin{equation} \label{tetracode}
\begin{array}{ll}
(1~3~2) (4~6~5)  (7~8~12), \qquad &  (1~3~2) (4~5~6)  (9~11~10), \\
(1~2~3) (7~8~12) (9~11~10), \qquad & (4~5~6) (7~8~12) (9~10~11),
\end{array}
\end{equation}
and their inverses.  They may be produced respectively by the paths
$$
\begin{array}{ll}
\left[\,0,12,1,0,9,6,11,10,5,0\,\right], \qquad &
\left[\,0,1,10,0,6,12,7,4,8,0 \,\right], \\
\left[\,0,12,1,0,9,5,6,11,4,0 \,\right], \qquad &
\left[\,0,12,10,0,3,4,2,1,5,0 \,\right],
\end{array}
$$
and their reverses.  Together with the identity,
these eight permutations form an elementary abelian group~$T$.
Notice that the orbits of the action of~$T$\/ on~$\PP \sm \{0\}$
are the sets $\ell \sm \{0\}$ for $\ell \in \starpt{0}$.
For any two distinct elements $\sigma,\tau \in T$,
there is exactly one line $\ell \in \LL(0)$
such that $\sigma(p)=\tau(p)$ for all $p\in\ell$.
Thus $T$\/ is the {\it tetracode\/}~\cite[p.~81]{Conway:1999}.
The elements of~$T$\/ are at maximal distance
not only from the identity but from each other (because $T$\/ is a group).

Since $d(\sigma) \leq 9$ for all $\sigma \in M_{13}$,
there is a much better algorithm than brute-force search
for ``solving'' the basic game---that is,
finding a short path producing a given permutation $\sigma \in M_{13}$.
It suffices to consider the case $\sigma \in M_{12}$,
since we can always start by moving the hole to $0$.

\begin{enumerate}
\item
Check if $\sigma \in N$, using the table (\ref{tetracode}).
If so, we are done.  If not, then $d(\sigma) \leq 8$.
\item
Create a list $L_1$ of all elements of $M_{13}$ of depth~$\leq 4$,
together with paths realizing them.
(An upper bound for the size of this list is $12 \cdot 9^3 = 8748$,
the number of paths $[0,p_1,\dots,p_4]$
with no three consecutive points collinear.)
\item
Create a list $L_2 = \{\sigma^{-1} \circ \tau \st \tau \in L_1\}$
of all elements of $M_{13}$ at distance~$\leq 4$ from $\sigma$.
\item
Since $d(\sigma) \leq 8$, we must have $L_1 \cap L_2 \neq 0$,
i.e., there are permutations $\tau$ and~$\tau'$
such that $\tau = \sigma^{-1}\tau'$.
Thus $\sigma = \tau' \tau^{-1}$,
and we can construct a path realizing $\sigma$
by concatenating those for $\tau'$ and $\tau^{-1}$.
\end{enumerate}

%%%%%%%%%%%%%%%%%%%%%%%%%%%%%%%%%%%%%%%%%%%%%%%%%%%%%%%%%%

\subsection{The signed game}

We now study the Cayley graph $\Cayley^+$ of the signed game,
whose vertices are the $2(13!/7!)$ positions of the signed game
and whose edges are given by signed moves.
As before, we can define distance, depth,
and numbers $\numd{2M_{12}}{k}$ and $\numd{2M_{13}}{k}$.

Let $\sigma \in M_{13}$, and let $\sigma_1,\sigma_2$
be the two lifts of $\sigma$ in~$2M_{13}$.
Then it is easy to see that
\begin{equation} \label{lift-depth}
d(\sigma) \ = \ \min(d(\sigma_1),d(\sigma_2)).
\end{equation}

\begin{prop} \label{depth-signed-prop}
The depth distributions for $2M_{12}$ and $2M_{13}$
are given by the following table:
$$
\begin{array}{l}
\begin{array}{c|ccccccc}
k & 0 & 1 & 2 & 3 & 4 & 5 & 6 \\ \hline
\numd{2M_{12}}{k} & 1 & 0 & 0 & 54 & 540 & 5184 & 25821 \\
\numd{2M_{13}}{k} & 1 & 12 & 108 & 918 & 7344 & 57852 & 356949 \\
\end{array}
\\ \\
\begin{array}{c|ccccccc}
k & 7 & 8 & 9 & 10 & 11 & 12 \\ \hline
\numd{2M_{12}}{k} & 85230 & 72351 & 898 & 0 & 0 & 1 \\
\numd{2M_{13}}{k} & 1192770 & 843291 & 11674 & 108 & 12 & 1
\end{array}
\end{array}
$$
\end{prop}

The unique element at maximal depth is $-\Id$,
the permutation that flips every counter
  in place
(see Proposition~\ref{sub-mathieu-prop}).
The subgroup $\{\Id,-\Id\}$ of $2M_{12}$ is central,
so every $\sigma \in 2M_{13}$ has a unique ``antipode''
$-\sigma = -\Id \cdot \sigma$,
which moves the counters to the same locations as~$\sigma$
but reverses all orientations,
and is uniquely maximally distant from~$\sigma$. 
Thus $\Cayley^+$ may be visualized as a ``globe''
in which pairs of poles represent antipodal permutations.

The depth distributions for $2M_{12}$ and $2M_{13}$
are the same as those for $M_{12}$ and $M_{13}$
for all depths $\leq 5$.
Indeed, let $k$ be the smallest number
such that $\numd{2M_{13}}{k} > \numd{M_{13}}{k}$.
By~(\ref{lift-depth}) and the pigeonhole principle,
there must be two elements $\sigma_1,\sigma_2 \in 2M_{13}$
at depth $\leq k$\/ which are lifts of the same $\sigma \in M_{13}$.
Then $\sigma_1^{-1}\sigma_2$ is a path of length $\leq 2k$\/
which induces the permutation $-\Id$ in some conjugate of $2M_{12}$,
which implies that $d(-\Id) \leq 2k$.
We must therefore have $k \geq 6$.

We also note that the depth distributions of $2M_{12}$ and $2M_{13}$
are ``symmetric near the poles'': there are the same numbers
of permutations at depths $0,1,2$ as at depths $12,11,10$ respectively.
However, the symmetry breaks down further from the poles:
fewer elements of $2M_{13}$ lie at depths $3,4,5$
than at depths $9,8,7$ respectively.
We may partially explain this phenomenon by noting that
\begin{equation} \label{depth-inequality}
d(\sigma) + d(-\sigma) \:\geq\: 12,
\end{equation}
for all $\sigma \in 2M_{13}$,
for otherwise $-\Id = -\sigma \circ \sigma^{-1}$ could be obtained
by a path of length strictly less than $12$.
Moreover, equality holds in (\ref{depth-inequality})
if and only if some minimal path from $\sigma$ to $-\sigma$
has $\Id$ as an intermediate position, which is not always the case.
Thus the mean depth of a permutation is greater than $6$.

Once again, these facts are based on the computational observation
that $-\Id$ is the unique element at depth~$12$.
This observation is also of use
in explaining the symmetry of the depth distribution near the poles.

\begin{prop}
Let $\sigma \in 2M_{13}$, with $d(\sigma) \in \{1,2\}$.
Then $d(-\sigma) = 12-d(\sigma)$.
\end{prop}

\begin{proof} Suppose that $d(\sigma)=1$;
then $\sigma$ is realized by a move sequence $[0,p]$, with $p \neq 0$.
Recall from (\ref{minus-id-path}) that $-\Id$
is realized by a \hbox{length-$12$} move sequence $[0,\dots,5,6,3,0]$.
Since $\Aut(\Prp)$ acts doubly transitively on~$\PP$,
we may choose $\alpha \in \Aut(\Prp)$
such that $\alpha(0)=0$ and $\alpha(3)=p$.  
Applying $\alpha$ to the move sequence realizing $-\Id$, we obtain
\begin{equation} \label{apath}
[0,\,\dots,\,\alpha(5),\,\alpha(6),\,p,\,0]
\end{equation}
which induces the signed permutation
$\alpha \circ -\Id \circ \alpha^{-1} = -\Id$.
Therefore, the path
$$
[0,\,\dots,\,\alpha(5),\,\alpha(6),\,\alpha(3)=p,\,0,\,p]
$$
induces the permutation $-\sigma$.
Deleting the last two moves, we obtain an equivalent path
of length~$11$.  So $d(\-\sigma) \leq 11$.
The opposite inequality follows from (\ref{depth-inequality}).

Similarly, if $d(\sigma)=2$,
then $\sigma$ is realized by a move sequence $[0,p,q]$,
with $0,p,q$ noncollinear.
By Proposition~\ref{P3-uniqueness}~(2),
$\Aut(\Prp)$ acts transitively on noncollinear triples of points,
so we may choose $\alpha \in \Aut(\Prp)$ such that
$\alpha(0)=0$, $\alpha(3)=p$, and $\alpha(6)=q$.
By the same argument as before,
$-\sigma$ is realized by the move sequence
$$
[0,\,\dots,\,\alpha(5),\,\alpha(6)=q],
$$
which has length~$10$.
\end{proof}

%%%%%%%%%%%%%%%%%%%%%%%%%%%%%%%%%%%%%%%%%%%%%%%%%%%%%%%%%%%%

\bibliographystyle{abstract}
  %% requires the file abstract.bst, which may or may not be
  %% a standard part of LaTeX distributions; available at
  %%     http://www.tex.ac.uk/tex-archive/biblio/bibtex/contrib/
\bibliography{biblio}
\end{document}